\documentclass[12pt]{amsart}

\usepackage{amssymb}
\usepackage{amsfonts}
\usepackage{amsmath}
\usepackage{latexsym}

\def\qed{\hspace{3.5mm} \hfill \vbox{\hrule height 3pt depth 2 pt width 2mm}
\bigskip}

\def\sfrac#1#2{{\scriptstyle\frac{#1}{#2}}}

\def\<{\langle}
\def\>{\rangle}

\def\Q{{\mathbb Q}}
\def\Z{{\mathbb Z}}\def\P{{\mathbb P}}

\def\goth{\mathfrak}
\def\ad{{\rm ad}}
\def\L{{\mathcal L}}
\def\F{{\mathcal F}}

\def\Der{{\rm Der}}

\def\drx{\partial_x}
\def\ddr#1{\frac{\partial}{\partial #1}}
\def\ddr#1{\partial_{#1}}

\def\Proof{\noindent \it Proof -- \rm}

\newtheorem{example}{Example}[section]

\newtheorem{corollary}[example]{Corollary}

\newtheorem{proposition}[example]{Proposition}

\title{Period polynomials and Ihara brackets}

\author[J.-G. Luque, J.-C. Novelli, and J.-Y. Thibon]%
{Jean-Gabriel Luque, Jean-Christophe Novelli,\\ and Jean-Yves Thibon}

\address{Institut Gaspard Monge, Universit\'e de
Marne-la-Vall\'ee \\
5, Boulevard Descartes \\Champs-sur-Marne \\77454 Marne-la-Vall\'ee cedex 2 \\
FRANCE}
\email[Jean-Gabriel Luque]{luque@univ-mlv.fr}
\email[Jean-Christophe Novelli]{novelli@univ-mlv.fr}
\email[Jean-Yves Thibon]{jyt@univ-mlv.fr}

\begin{document}

\begin{abstract}
Schneps [J. Lie Theory {\bf 16} (2006), 19--37] has found surprising links
between Ihara brackets and even period polynomials. These results can be
recovered and generalized by considering some identities relating Ihara
brackets and classical Lie brackets. The period polynomials generated by this
method are found to be essentially the Kohnen-Zagier polynomials.
\end{abstract}

\maketitle

\section{Introduction}

The Ihara bracket is defined on a subspace of the free Lie algebra $\L$
on two generators, whose elements are intepreted as special derivations of
$\L$.
One can find relations between the two brackets. This gives rise to
interesting identities, which surprisingly, turn out to be related to period
polynomials~\cite{Sch}.

The \emph{period polynomials} $r^\pm(f)$ of a modular form
$f\in M_k(SL_2(\Z))$ are defined by
\begin{equation}
r(f)(t)=\int_0^{i\infty}f(z)(t-z)^w dz \quad (w=k-2)\,,
\end{equation}
and $r=r^+ + ir^-$. These polynomials satisfy the \emph{period relations}
\begin{equation}
P(t)+t^wP\left(\frac{-1}{t}\right) = 0\,,
\end{equation}
\begin{equation}
P(t)+t^wP\left(1-\frac{1}{t}\right)+(t-1)^wP\left(\frac{1}{1-t}\right)=0\,.
\end{equation}
The solutions of this system are in one-to-one correspondence with cusp
forms (except for the $t^w-1$ which are obtained from Eisenstein series).
This is the Eichler-Shimura correspondence (see, \emph{e.g.}, \cite{KZ}).

Writing
\begin{equation}
r(f)(t)=\sum_{n=0}^w i^{-n+1}\binom{w}{n}r_n(f)\,,
\end{equation}
one defines linear forms $r_n$ on $M_k$. Since this is an inner product space
for the Petersson scalar product, there exist modular forms $R_{n;k}$ such
that for all $f\in S_k$,
\begin{equation}
r_n(f)=(f,R_{n;k})\,.
\end{equation}
The period polynomials of the $R_{n;k}$ have been explicitly computed
by Kohnen and Zagier \cite[Theorem 1, p. 208]{KZ}. Their result
can be stated as follows:
let $B_n(t)$ be the Bernoulli polynomials and $b_n$ be the Bernoulli numbers
defined by
\begin{equation}
\frac{u e^{tu}}{e^u-1} =: \sum_{n\geq0} B_n(t) \frac{u^n}{n!}
\text{\qquad and\qquad}
B_n(t) =: \sum_{i=0}^n \binom{n}{i} b_i t^{n-i}.
\end{equation}
and let
\begin{equation}
\begin{split}
P^\pm_{n;k}(t) :=& \pm\frac{1}{n+1}\left[B_{n+1}(t)\mp
t^wB_{n+1}\left(\frac1t\right)\right]
\\&+\frac{1}{p+1}\left[B_{p+1}(t)\mp
t^wB_{p+1}\left(\frac1t\right)\right]\,,
\end{split}
\end{equation}
where $w=k-2=n+p$. Then \cite{KZ},
\begin{equation}
P^\pm_{n;k}(t)= r^\pm(\lambda^\pm_{n;k} R_{n;k} +\mu^\pm_{n;k}G_k)
\end{equation}
where $G_k$ is the Eisenstein series and $\lambda^\pm_{n;k},\mu^\pm_{n;k}$
are explicit constants that we shall not need. 

Hence, any linear combination of the $P^\pm_{n;k}(t)$ for a given $k$ is a
period polynomial.
Such combinations arise naturally from the Ihara bracket. In
\cite{Sch}, Schneps obtained a characterization of even period polynomials in
terms of linear relations between certain brackets.
In this note, we give a new proof of this result and obtain a similar
characterization of odd period polynomials. We also give a simpler way
to generate these polynomials.

\medskip
{\it Acknowledgements.-}
The authors are grateful to Don Zagier for recognizing their period
polynomials as those of~\cite{KZ}.

\section{Ihara brackets and the stable derivation algebra}
\label{stable}

The stable derivation algebra ${\goth F}$ appears in works by Ihara on Galois
representations on $\pi_1(\P^1_\Q-\{0,1,\infty\})$
\cite{Ih1,Ih2,Ih3,Ih4,Ih5,Ih6}. It plays also a crucial r\^ole in recent works
on multiple zeta values \cite{Fur,Gon,Rac}.

The underlying vector space of ${\goth F}$ can be identified with a
subalgebra of the free Lie algebra over two letters \cite{Ih6}.
Let $\L$ be the free $\Q$-Lie algebra  over two letters
$a,b$. By Lazard elimination \cite{Reut}, the Lie algebra decomposes as 
\begin{equation}
\L=\Q a\oplus \F^1\L
\end{equation}
where the Lie algebra $\F^1\L$ is free over the infinite sequence 
\begin{equation}
\phi_{n+1}=\frac1{n!} \ad^n_a(b)\qquad  (n\geq0),
\end{equation}
where $\ad$ denotes the adjoint representation.

Let $V$ be any vector space and $\Der(\L(V))$ be the Lie algebra of
derivations of the free Lie algebra $\L(V)$.
The commutator of two inner derivations satisfies
\begin{equation}
[\ad_f,\ad_g] = \ad_{[f,g]}
\end{equation}
where $f,g\in\L$ and $[f,g]$ is the bracket in $\L$.
More generally, if $\pi:\ V\rightarrow V$ is a projector, one can associate
with any element $f\in\L(V)$ a derivation $D_f$ defined by
\begin{equation}
D_f(v)=[\pi(v),f]\,,
\end{equation}
which satisfies
\begin{equation}
[D_f,D_g] = D_{\{f,g\} }\,
\text{\ \ where\ \ }
\{f,g\} := [f,g] + D_g(f) - D_f(g)\,.
\end{equation}
In the special case where $V=\Q a\oplus \Q b$ and $\pi(a)=a$, $\pi(b)=0$, the
operation $\{f,g\}$ is called the Ihara bracket.
Let $\mathcal G$ be the Lie subalgebra of $\Der(\L)$ generated by the special
derivations $D_f$, $f\in\L$.
It can be identified with the Ihara subalgebra $\F^1\L$ of $\L$ generated by
the $\phi_n$, so that it is also a Lie subalgebra of $\L$ for the ordinary
bracket.
The stable derivation algebra $\goth F$ is a subalgebra of $\mathcal G$,
spanned by the $D_f$ such that $f$ satisfies certain
identities~(\cite{Ih6}, see also~\cite{Fur}).

In \cite{Fur,Ih4}, Furusho and Ihara consider the filtration 
\begin{equation}
\L = \F^0 \L \supset\F^1 \L\supset\cdots\supset\F^n \L\supset\cdots
\end{equation}
where 
\begin{equation}
\F^0 \L = \L\ \text{and $\F^n \L = [\F^1 \L, \F^{n-1} \L]$ for $n\geq 2$}.
\end{equation} 
This induces a filtration on ${\goth F}$ 
\begin{equation}
{\goth F} = \F^1{\goth F}\supseteq \F^2{\goth F}
\supseteq\cdots\supseteq \F^n{\goth F}\supseteq \cdots\,, 
\end{equation}
and the compatibility of the Ihara bracket with this filtration implies that
the Lie algebra $\goth F$ is generated by polynomials which do not belong to
$\F^2{\goth F}$.
From \cite{Ih6}, 
$\dim_\Q \F^1{\goth F}_n/ \F^2{\goth F}_n=1$ if $n$ is odd and
greater than $1$ and $0$ otherwise. 
Hence, the Lie algebra $\goth F$ is generated by a set
\begin{equation}
\{f_{2n+1}\in {\goth F}_{2n+1}\,|\,n\geq 1\}.
\end{equation} 
such that
\begin{equation}\label{modF2}
f_{2n+1}\equiv (2n)!\,\phi_{2n+1}\mbox{ mod } \F^2{\goth F}
\end{equation}
Therefore,
\begin{equation}\label{modF3}
\{f_{2n+1},f_{2n'+1}\} \equiv (2n)!\,(2n')!\,\{\phi_{2n+1},\phi_{2n'+1}\}
\mbox{ \ mod } \F^3{\goth F}
\end{equation}
which implies the equivalence between equalities of the form 
\begin{equation}
\sum_{\genfrac{}{}{0pt}{}{i+j=n}{i,j\leq 3\text{ odd}}}a_{ij}\,\{f_i,f_j\}
\equiv 0\mbox{ mod } \F^3{\goth F}
\end{equation}
and 
\begin{equation}\label{aFF}
\sum_{\genfrac{}{}{0pt}{}{i+j=n}{i,j\leq 3\text{ odd}}}a_{ij}\,(i-1)!\,(j-1)!\,
\{\phi_i,\phi_j\}= 0\,. 
\end{equation}
This is precisely the kind of linear relations giving rise to period
polynomials. Those found by Schneps \cite{Sch} are obtained by substituting
\begin{equation}\label{t-t}
\{\phi_i,\phi_j\} \rightarrow \frac{1}{(i-1)!(j-1)!}(t^{i-1}-t^{j-1})
\end{equation}
in the left-hand side of relations like (\ref{aFF}).
Let us note that if we set
$\varphi_i(t)=\frac1{(i-1)!}t^{i-1}$ and 
$\<f(t),g(t)\>=f(t)g(1)-f(1)g(t)$ (this is a Lie bracket),
then (\ref{t-t}) amounts to replace $\{\phi_i,\phi_j\}$
by $\<\varphi_i(t),\varphi_j(t)\>$.

All the calculations presented in this paper aim at finding
explicit linear relations between the 
$[\phi_m,\phi_n]$, $\{\phi_m,\phi_n\}$ and $D_{\phi_m}(\phi_n)$.
As a general rule, the coefficients are expressed in terms of Bernoulli
numbers and can be related to the Kohnen-Zagier polynomials.

\section{Computation of Ihara brackets}
\label{Ihara}

The aim of this section is to express the Ihara bracket in terms of
classical brackets. Our first step is to compute the action of the special
derivations on the generators $\phi_n$.
Let us start with the simple equality
\begin{equation}\label{Dptop}
D_{\phi_n}a=-n\phi_{n+1}\,.
\end{equation}
Since $D_{\phi_n}$ is a derivation, it follows from (\ref{Dptop}) that
\begin{equation}
\label{expand1}
D_{\phi_n}(\phi_p) =
  -\frac{n}{(p-1)!}
  \sum_{i=1}^{p-1}(p-i-1)!\,\ad_a^{i-1}\,[\phi_{n+1},\phi_{p-i}].
\end{equation}
The Leibniz formula gives
\begin{equation}
\ad_a^{i-1}[\phi_{n+1},\phi_{p-i}]
=\displaystyle\sum_{j=0}^{i-1}\binom{i-1}{j}
 \frac{(n+j)!(p-j-2)!}{n!(p-i-1)!}\,[\phi_{n+j+1},\phi_{p-j-1}]\,.
\end{equation}
Substituting  in (\ref{expand1}) and rearranging the sums, one gets, for all
$n,p\geq1$
\begin{equation}
\begin{split}
\label{Dptop2}
D_{\phi_n}(\phi_p)
&=\sum_{j=1}^{p-1}\binom{n+j-1}{j}[\phi_{p-j},\phi_{n+j}]\\
&=\sum_{k=1}^{p-1}\binom{n+p-1-k}{p-k}[\phi_k,\phi_{p+n-k}]\, .
\end{split}
\end{equation}
Hence, from the definition of Ihara bracket, one has, for all $n,p\geq1$
\begin{equation}
\begin{split}
\{ \phi_n,\phi_p \}
=&\ \ [\phi_n,\phi_p] \\
& + \sum_{k=1}^{n-1} \binom{n+p-1-k}{n-k} [\phi_k,\phi_{n+p-k}] \\
& - \sum_{k=1}^{p-1} \binom{n+p-1-k}{p-k} [\phi_k,\phi_{n+p-k}]\,. \\
\end{split}
\end{equation}
This equation can be rewritten as
\begin{equation}
\label{PtoF}
\{ \phi_n,\phi_p \} =
\sum_{k=1}^{\max(n\!-\!1,p\!-\!1)}
\left( \binom{n\!+\!p\!-\!1\!-\!k}{n-k} - \binom{n\!+\!p\!-\!1\!-\!k}{p-k}
\right)
[\phi_k,\phi_{n+p-k}]\,,
\end{equation}
for all $n,p\geq1$ as one can check on all cases $n<p$, $n=p$, and $n>p$.

In terms of generating series, both Equations~(\ref{Dptop2}) and~(\ref{PtoF})
have simple expressions. Set
\begin{equation}\label{genp}
\Phi(x):=\sum_{n\geq 1}\phi_nx^{n-1}\,.
\end{equation}
Then,
\begin{equation}
\begin{split}
D_{\Phi(x)}\Phi(y) &:=
    \sum_{n,p\geq 1}D_{\phi_n}(\phi_p)\, x^{n-1}y^{p-1}\\
 &=
    \sum_{n,p\geq 1} \sum_{k=1}^{p-1}\binom{n+p-1-k}{p-k}
                    [\phi_{k}\,y^{k-1},\phi_{n+p-k}\,x^{n-1}y^{p-k}]\\
 &=\left[
         \sum_{k\geq 1}\phi_k\, y^{k-1},\ \
          \sum_{s\geq 1}\phi_s\sum_{n=1}^{s-1}\binom{s-1}{n-1}
               x^{n-1}y^{s-n}\right]\\
 &= \left[\Phi(y),\sum_{s\geq1}\phi_s((x+y)^{s-1}-x^{s-1})\right]\\
 &= [\Phi(y),\Phi(x+y)-\Phi(x)]\\[10pt]
 &= [\Phi(x),\Phi(y)]+[\Phi(y),\Phi(x+y)].
\end{split}
\end{equation}
Hence,
\begin{equation}\label{sdphi}
D_{\Phi(x)}\Phi(y)=[\Phi(x),\Phi(y)]+[\Phi(y),\Phi(x+y)]
\end{equation}
and
\begin{equation}\label{ihlie}
\{\Phi(x),\Phi(y)\}=[\Phi(y),\Phi(x)]+[\Phi(x)-\Phi(y),\Phi(x+y))]\,.
\end{equation}

\section{Inversion of Equation~(\ref{sdphi}) and period polynomials}

\subsection{Generic inversion of~(\ref{sdphi})}

Let $(F_{i,j})_{i,j\geq1}$ and $(G_{i,j})_{i,j\geq1}$ be two bi-indexed
sequences of elements of some vector space whose generating series
\begin{equation}
F(x,y) = \sum_{i,\,j\geq1} F_{i,j}\, x^{i-1} y^{j-1},
\qquad
G(x,y) = \sum_{i,\,j\geq1} G_{i,j}\, x^{i-1} y^{j-1},
\end{equation}
satisfy
\begin{equation}
\label{F2G}
F(x,y) = G(x,y) - G(x+y,y) = (1-e^{y\ddr{x}})\, G(x,y).
\end{equation}
Thanks to Equation~(\ref{sdphi}), this is the case of
\begin{equation}
F(x,y)=D_{\Phi(x)}\Phi(y) \text{\quad and\quad }
G(x,y)=[\Phi(x),\Phi(y)].
\end{equation}
%
If $F(x,y)$ is given, this formula determines $G(x,y)$ up to a function of
$y$:
\begin{equation}
\begin{split}
y\drx G(x,y) =& \frac{y\drx}{1-e^{y\drx}} F(x,y)\\
&= - \sum_{k\geq0} \frac{b_k}{k!} y^k \drx^k
\sum_{i,j\geq1}F_{i,j}\, x^{i-1}y^{j-1} \\
&= - \sum_{k\geq0} \frac{b_k}{k!}
  \sum_{i\geq k+1;\,j\geq1} (i-1)\cdots(i-k) F_{i,j}\, x^{i-1-k} y^{j-1+k} \\
&= - \sum_{k\geq0}
  \sum_{i\geq k+1;\,j\geq1} \binom{i-1}{k}b_k\, F_{i,j}\,
x^{i-1-k} y^{j-1+k} \\
&= - \sum_{i\geq1} \sum_{k=0}^{i-1} \binom{i-1}{k}b_k\, x^{i-1-k}y^k
\sum_{j\geq1} F_{i,j}\,y^{j-1} \\
&= - \sum_{i\geq0} \sum_{k=0}^{i} \binom{i}{k}b_k\,
\left(\frac{x}{y}\right)^{i-k} y^i \sum_{j\geq1} F_{i+1,j}\,y^{j-1}
\end{split}
\end{equation}
so that, finally
\begin{equation}
\label{FG}
y\drx G(x,y) = - \sum_{i\geq0} B_i\left(\frac{x}{y}\right) y^i 
\sum_{j\geq1} F_{i+1,j}\,y^{j-1}.
\end{equation}

Then, comparing the coefficients of $x^ny^p$ (with $n\geq0$ and $p\geq1$) on
both sides of Equation~(\ref{FG}), one obtains

\begin{equation}
\begin{split}
(n+1) G_{n+2,p}
&= - \sum_{i=n}^{n+p} \binom{i}{i-n} b_{i-n}\, F_{i+1,n+p+1-i} \\
&= - \sum_{i=0}^p     \binom{n+i}{i} b_i    \, F_{n+1+i,p+1-i}\,,
\end{split}
\end{equation}
so that

\begin{proposition}
For all $n\geq 2$ and $p\geq1$,
\begin{equation}
\label{GF}
G_{n,p} = - \frac{1}{n-1}
\sum_{i=0}^p \binom{n-2+i}{n-2} b_i \, F_{n-1+i,p+1-i}\,.
\end{equation}
\qed
\end{proposition}

Now, assume that $G_{n,p}=[\phi_n,\phi_p]$, so that
$F_{n,p}=D_{\phi_n}(\phi_p)$.

\begin{proposition}
For all $n\geq 2$ and $p\geq1$,
\begin{equation}
\label{Dpcroch}
\begin{split}
[\phi_n,\phi_p]
=& - \frac{1}{n-1}
\sum_{i=0}^{p}\binom{n-2+i}{n-2} b_i \, D_{\phi_{n-1+i}}(\phi_{p+1-i}) \\
=& - \sum_{i=0}^p \binom{n-1+i}{n-1} \frac{b_i}{n\!-\!1\!+\!i}
D_{\phi_{n-1+i}}(\phi_{p+1-i}).
\end{split}
\end{equation}
\qed
\end{proposition}

Note that in the first equation, the summation can be taken up to
$p-1$ since $D_{\phi_i}(\phi_1)=D_{\phi_i}(b)=0$ for all $i$.
Now, since the Lie bracket is antisymmetric, one gets the following
relations between the $D_{\phi_i}(\phi_j)$:
\begin{corollary}
\label{cor-annul-dphi}
~

\noindent
(i) For any $n>1$, 
\begin{equation}
\label{Dpcroch01}
\sum_{i=0}^{n}\binom{n-2+i}{n-2} b_i\, D_{\phi_{n-1+i}}(\phi_{n+1-i}) = 0,
\end{equation}
(ii) Setting $b_i=0$ if $i<0$, one has, for all $n,p\geq2$
\begin{equation}
\label{Dpcroch02}
\sum_{i=1}^{n+p-1}
\left( \binom{i-1}{i-p+1}\frac{b_{i-p+1}}{p-1}
     + \binom{i-1}{i-n+1}\frac{b_{i-n+1}}{n-1}\right)
D_{\phi_i}(\phi_{n+p-i})=0.
\end{equation}
\qed
\end{corollary}

\subsection{Period polynomials}
\label{period1}

Let us now consider the specialization
\begin{equation}
F_{i,j} = \frac{1}{(i-1)!} \frac{1}{(j-1)!} (t^{i-1} - \epsilon t^{j-1}).
\end{equation}
Then, Equation~(\ref{GF}) gives
\begin{equation}
\begin{split}
G_{n,p} =& \frac{-1}{n-1}
\sum_{i=0}^{p}\binom{n-2+i}{n-2} b_i \,
\frac{t^{n-2+i}-\epsilon t^{p-i}}{(n-2-i)!(p-i)!} \\
=& \frac{-1}{(n-1)!p!}
   \sum_{i=0}^p \binom{p}{i} b_i (t^{n-2+i}-\epsilon t^{p-i})\\
=& \frac{-1}{(n-1)!p!}
   \left(t^{n+p-2}B_p\left(\frac1t\right) - \epsilon B_p(t)\right)\\
=& \frac{1}{(n-1)!(p-1)!}
     \left(\frac1p \left(\epsilon B_p(t) - t^{n+p-2}B_p\left(\frac1t\right)
\right)\right).\\
\end{split}
\end{equation}
Note that in terms of generating series, the expressions of $F$ and $G$ are
simple:
\begin{equation}
F(x,y) = e^{tx+y} - \epsilon e^{x+ty}
\text{ and }
G(x,y) = f(y,t) + \frac{e^{tx+y}}{1-e^{ty}} - \epsilon \frac{e^{x+ty}}{1-e^y}.
\end{equation}

One recognizes in the coefficients $G_{n,p}$ the building blocks of the period
polynomials introduced by Kohnen and Zagier~(\cite{KZ}, Theorem 1). It
follows from their results that, with $\epsilon =1$,
\begin{equation}
n!\,p!\, (G_{n+1,p+1} + G_{p+1,n+1}) = P^+_{n,n+p+2}(t)
\end{equation}
is an even period polynomial for $n$, $p$ even, and with $\epsilon=-1$,
\begin{equation}
n!\,p!\, (G_{n+1,p+1} - G_{p+1,n+1}) = P^-_{n,n+p+2}(t)
\end{equation}
is an odd period polynomial for $n$, $p$ odd.

Note that in the case of even period polynomials, this amounts to substituting
\begin{equation}
D_{\phi_i}(\phi_j) \mapsto \frac{1}{(i-1)!(j-1)!} (t^{i-1}-t^{j-1})
\end{equation}
in the left-hand sides of the linear relations~(\ref{Dpcroch01})
and~(\ref{Dpcroch02}), which is analogous to the result of Schneps~\cite{Sch}.

\section{Ordinary brackets in terms of Ihara brackets}
\label{inv2}

We shall now give a partial inversion of Equation~(\ref{PtoF}), \emph{i.e.},
express the Lie brackets $[\phi_{2n},\phi_k]$ as a linear combination of Ihara
brackets. Here are some examples:
\begin{equation}
[\phi_2,\phi_k]=
  \sfrac{1}{k-1} \{\phi_3,\phi_{k-1}\} - \sfrac{1}{2} \{\phi_2,\phi_k\}.
\end{equation}
\begin{equation}
[\phi_4,\phi_k]=
  \sfrac{k}{12}  \{\phi_3,\phi_{k+1}\}
- \sfrac{1}{2}   \{\phi_4,\phi_k\}
+ \sfrac{1}{k-1} \{\phi_5,\phi_{k-1}\}.
\end{equation}
\begin{equation}
\begin{split}
[\phi_6,\phi_k]= &
- \sfrac{(k+2)(k+1)k}{720}\{\phi_3,\phi_{k+3}\}
+ \sfrac k{12}\{\phi_5,\phi_{k+1}\}
- \sfrac12\{\phi_6,\phi_k\} \\
&+ \sfrac1{k-1}\{\phi_7,\phi_{k-1}\}.
\end{split}
\end{equation}
\begin{equation}
\begin{split}
[\phi_8,\phi_k] =&
  \sfrac{(k+4)(k+3)(k+2)(k+1)k}{30240}\{\phi_3,\phi_{k+5}\}
- \sfrac{(k+2)(k+1)k}{720}\{\phi_5,\phi_{k+3}\}\\ 
& + \sfrac k{12}\{\phi_7,\phi_{k+1}\}
- \sfrac12\{\phi_8,\phi_k\}
+ \sfrac1{k-1}\{\phi_9,\phi_{k-1}\}.
\end{split}
\end{equation}
The general formula is as follows.
\begin{proposition}
For each $k\geq2$ and $n\geq1$, one has
\begin{equation}
\label{LToP}
[\phi_{2n},\phi_k]
=\sum_{i=0}^{2n} \binom{k-1+i}{k-1} \frac{b_i}{k-1+i}
  \{\phi_{2n-i+1},\phi_{k+i-1}\}.
\end{equation}
\end{proposition}
\Proof
The generating series of the right-hand side
of (\ref{LToP}) is
\begin{equation}
\label{dS}
S(x,y)=\sum_{n\geq 1}\sum_{k\geq 2} x^{2n-1}y^{k-1}
        \sum_{i=0}^{2n-1}\frac{(k+i-2)!}{(k-1)!}\frac{b_i}{i!}
                         \{\phi_{2n-i+1},\phi_{k+i-1}\}.
\end{equation}
%
Rearranging the sum, one obtains
\begin{equation}
\begin{split}
x\ddr{y} S(x,y)=&
  \sum_{n\geq1}\sum_{k\geq2} x^{2n}(k\!-\!1)y^{k-2}
    \sum_{i=0}^{2n} \frac{(k+i-2)!}{(k-1)!}\frac{b_i}{i!}
                    \{\phi_{2n-i+1},\phi_{k+i-1}\} \\
=& \sum_{n\geq1}\sum_{k\geq2} \sum_{i=0}^{2n}
   \frac{(k+i-2)!}{(k-2)!}\frac{b_i}{i!}
   \{\phi_{2n-i+1}x^{2n},\phi_{k+i-1}y^{k-2}\}\\
=& \sum_{i\geq0} x^i \frac{b_i}{i!}
   \left\{ \sum_{n\geq\lceil i/2\rceil}\phi_{2n\!-\!i\!+\!1}x^{2n-i},
      \sum_{k\geq2} \frac{(k\!+\!i\!-\!2)!}{(k-2)!}y^{k-2} \phi_{k+i-1}\right\}\\
=& \sum_{i\geq0} x^i \frac{b_i}{i!}
   \left\{ \frac{1}{2}(\Phi(x)+ (-1)^i\Phi(-x)),
      \sum_{k\geq2} \ddr{y}^i\, y^{k+i-2} \phi_{k+i-1}\right\}\\
=& \frac{1}{2}
\sum_{i\geq 0}\left(x\ddr{y}\right)^i
\frac{b_i}{i!}\{\Phi(x)+(-1)^i\Phi(-x),\Phi(y)\}\,,
\end{split}
\end{equation}
so that
\begin{equation}
\label{eqS1}
2x\ddr{y} S(x,y)=
 \frac{x\ddr{y}}{e^{x \ddr{y}}-1}\{\Phi(x),\Phi(y)\}
+ \frac{-x\ddr{y}}{e^{-x\ddr{y}}-1}\{\Phi(-x),\Phi(y)\}.
\end{equation}
Equation (\ref{ihlie}) gives
\begin{equation}
\{\Phi(x),\Phi(y)\}=
\left( e^{x\ddr{y}} + e^{y\ddr{x}} -1 \right)
[\Phi(x),\Phi(y)]\,.
\end{equation}
Substituting this expression in (\ref{eqS1}), one gets
\begin{equation}
\label{eqS2}
\begin{split}
2x\ddr{y} S(x,y)=&
\  x\ddr{y} \left(1+\frac{1}{e^{x\ddr{y}} -1} e^{y\ddr{x}} \right)
  [\Phi(x),\Phi(y)]\\
&
  -x\ddr{y}\left(1+\frac{1}{e^{-x\ddr{y}} -1} e^{-y\ddr{x}} \right)
  [\Phi(-x),\Phi(y)].\\
=&\ x\ddr{y}\Big([\Phi(x)-\Phi(-x),\Phi(y)]
         + \frac{1}{e^{x\ddr{y}}-1}  [\Phi(x+y),\Phi(y)] \\
 &\qquad - \frac{1}{e^{-x\ddr{y}}-1} [\Phi(-x+y),\Phi(y)]\Big)\\
=&\ x\ddr{y}\Big([\Phi(x)-\Phi(-x),\Phi(y)]
         + \frac{1}{e^{x\ddr{y}}-1} ( [\Phi(x+y),\Phi(y)]\\
 &\qquad + e^{x\ddr{y}}
    [\Phi(-x+y),\Phi(y)]\Big)\\[10pt]
=&\ x\ddr{y}\left([\Phi(x)-\Phi(-x),\Phi(y)]\right).\\
\end{split}
\end{equation}
Hence,
\begin{equation}
\begin{split}
\label{eqS3}
x\ddr{y} S(x,y) =& \frac12\, x\ddr{y}\, [\Phi(x)-\Phi(-x),\Phi(y)]\\
=&\sum_{n\geq1}\sum_{k\geq2}(k-1)\,[\phi_{2n},\phi_k]\,x^{2n}y^{k-2}.
\end{split}
\end{equation}
Comparing the coefficients of $x^{2n-1}y^{k-1}$ in (\ref{dS}) and of
$x^{2n}y^{k-2}$ in (\ref{eqS3}) for $n\geq1$ and $k\geq2$, one obtains
(\ref{LToP}).
\qed

Note that Formula~(\ref{LToP}) is very similar to Formula~(\ref{Dpcroch}) when
substituting $n=k$ and $p=2n$:
\begin{equation}
[\phi_{2n},\phi_k]
=\sum_{i=0}^{2n} \binom{k-1+i}{k-1} \frac{b_i}{k\!-\!1\!+\!i}
  \{\phi_{2n-i+1},\phi_{k+i-1}\}.
\end{equation}
\begin{equation}
\label{repri}
\begin{split}
[\phi_{2n},\phi_k]
= \sum_{i=0}^{2n} \binom{k-1+i}{k-1} \frac{b_i}{k\!-\!1\!+\!i}
D_{\phi_{k-1+i}}(\phi_{2n+1-i}).
\end{split}
\end{equation}
%
%
Thus, we obtain without further calculations the following analogs
of~(\ref{Dpcroch01}) and~(\ref{Dpcroch02}):

\begin{corollary}
For all $n\geq1$,
\begin{equation}
\label{cor1}
\sum_{i=0}^{2n}\binom{2n-2+i}{2n-2} b_i
\{\phi_{2n-i+1},\phi_{2n+i-1}\} = 0.
\end{equation}
For all $n\geq1$ and $p\geq1$,
\begin{equation}
\label{cor2}
\sum_{i=1}^{2n+2p-1}
\left( \binom{i-1}{i\!-\!2p\!+\!1}\frac{b_{i\!-\!2p\!+\!1}}{2p-1}
     + \binom{i-1}{i\!-\!2n\!+\!1}\frac{b_{i\!-\!2n\!+\!1}}{2n-1}\right)
\{ \phi_{2n+2p-i}, \phi_i \} =0.
\end{equation}
\end{corollary}

Note that the left-hand sides of Equations~(\ref{cor1}) and~(\ref{cor2}) are
only composed of brackets of odd $\phi$. For example, the first equation gives

\begin{equation}
\begin{split}
0=&\ 9\, \{\phi_5,\phi_7\} - 14 \{\phi_3,\phi_9\},\\
0=&\ 11\, \{\phi_7,\phi_9\} - 21\, \{\phi_5,\phi_{11}\}
    + 66\, \{\phi_3,\phi_{13}\},\\
0=&\ 13\, \{\phi_9,\phi_{11}\} - 33\,\{\phi_7,\phi_{13}\}
    + 143\, \{\phi_5,\phi_{15}\} - 858\, \{\phi_3,\phi_{17}\},\\
0=&\ 300\, \{\phi_{11},\phi_{13}\} - 1001\, \{\phi_9,\phi_{15}\}
    + 5720\, \{\phi_7,\phi_{17}\} - 43758\, \{\phi_5,\phi_{19}\}
   \\ &\  + 419900\, \{\phi_3 ,\phi_{21}\}.
\end{split}
\end{equation}
whereas the second one gives
\begin{equation}
\begin{split}
0=&\ 195\, \{\phi_{11},\phi_7\} - 825\, \{\phi_{13},\phi_5\}
   - 4004\,\{\phi_{15}, \phi_3\},\\
0=&\ 85\, \{\phi_{13},\phi_9\} - 442\, \{\phi_{15}, \phi_7\}
   + 2730\, \{\phi_{17}, \phi_5\} - 21216\, \{\phi_{19},\phi_3\},\\
0=&\ 2193\, \{\phi_{13}, \phi_{11}\} - 7973\, \{\phi_{15},\phi_9\}
   + 47213\, \{\phi_{17},\phi_7\}\\
  &\ - 364803\, \{\phi_{19}, \phi_5\}
   + 3509718\, \{\phi_{21}, \phi_3\}.
\end{split}
\end{equation}

\subsection{Period polynomials}
\label{period2}

\noindent\hskip.5cm
It follows from the discussion of Section~\ref{period1} that if one
writes~(\ref{cor1}) and~(\ref{cor2}) as
\begin{equation}
\label{Schn1}
\sum_{i,j\text{ odd}} a_{i,j} \{\phi_i,\phi_j\} =0\,,
\end{equation}
then
\begin{equation}
\label{Schn2}
\sum_{i,j\text{ odd}} a_{i,j} \frac{1}{(i-1)!(j-1)!}
                     \left( t^{i-1} - t^{j-1}\right)
\end{equation}
is a period polynomial as first shown by Schneps~\cite{Sch}. Actually, Schneps
has shown that a relation of the type~(\ref{Schn1}) holds iff~(\ref{Schn2}) is
a period polynomial.
Let us recall the explanation: Ihara and Takao~\cite{Ih6} have proved that the
space of linear relations of the form~(\ref{Schn2}) has the same dimension as
the space of cusp forms $S_n(SL_2(\Z))$.
We have seen that all the even period polynomials $P^+_{2n,2m}(t)$ can be
obtained in this way. Hence, all linear relations between the
$\{\phi_i,\phi_j\}$ (with $i,j$ odd) are consequences of~(\ref{Schn1}), so
that we recover the result of Schneps.


\end{document}